\documentclass[12pt,fleqn]{amsart}

\usepackage{latexsym,amssymb,amsmath,amsthm,amsaddr}
\usepackage{amscd}
\usepackage{enumerate}
\usepackage{amsfonts}
 \usepackage{anysize}
 \usepackage{bbm}
 \usepackage{graphicx}
 \usepackage{url}
 \usepackage{float}

\newtheorem{theorem}{Theorem}[section]
\newtheorem{lemma}[theorem]{Lemma}

\newtheorem{proposition}[theorem]{Proposition}

\theoremstyle{definition}
\newtheorem{definition}[theorem]{Definition}
\newtheorem{remark}[theorem]{Remark}
\newtheorem{example}[theorem]{Example}

\numberwithin{equation}{section}
\newcommand{\N}{{\mathbb N}}

\newcommand{\R}{{\mathbb R}}

\newcommand{\vr}{\varepsilon}

\newcommand{\fn}{\mathbb{E}\sp{1}}

\title[Neural network operators on fuzzy continuous functions]{Neural Network Operators on fuzzy number valued continuous functions}

\author{Juan J. Font, Sergio Macario, Manuel Sanchis}

 \address{Institut Universitari de Matem\`{a}tiques i Aplicacions de
Castell\'{o}
 (IMAC),\\ Universitat Jaume I, Campus del Riu Sec. s/n, 12071 Castell\'{o}
(Spain)}

 \email{font@uji.es, macario@uji.es, sanchis@uji.es}

 \subjclass[2010]{}

\newcommand{\restringido}[1]{\,\vrule height5pt width.4pt depth10pt\,\lower
10pt\hbox{\scriptsize $#1$}}

\begin{document}

\begin{abstract}
 We extend Cardaliaguet-Euvrard neural network operators to the context of fuzzy number valued continuous functions and study their behaviour. We focus on level continuous, sendograph continuous and endograph continuous functions and obtain Jackson-type results in all these cases.
 
 \end{abstract}

\maketitle

\section{Introduction}
A multilayer feedforward  artificial neural network can be described, in a simple way, by a mathematical function as
$$N(x)=\sum_{j=0}^{n} c_j\sigma(a_j\cdot x+b_j),$$
where $x\in \R^m$, $m\in\N$; $c_j,a_j,b_j\in\R$ are the coefficients, weights and threshold values, respectively, and $\sigma$, which is usually a sigmoidal function, is called the activation function.

In the last forty years, such neural networks have been successfully applied to approximate continuous functions of one or several variables, that is, they are universal approximators (\cite{cybenko, hornik}).

Neural network (NN) operators, which are an extension of artificial neural networks, were introduced in order to provide constructive approximation algorithms, making them more intuitive than
those used in the context of neural networks in the sense that all the essential components
of a neural network are known in these operators, e.g., coefficients, weights
and threshold values. These operators admit a function as an input and gives also a function as an output. The simplest way is given by
$$N_{n}(f,x)=\sum_{j=0}^n \phi_j(x)f\left(x_j\right)$$
where $f\in C([a,b])$,  $\{x_j\}_{j=0}^n$ a partition of $[a,b]$ and the coefficients $\phi_j(x)$ depend on the activation function.

The NN network operators we are interested in were introduced by Anastassiou (\cite{anas}) based on the ideas proposed by Cardaliaguet and Euvrard in \cite{carda}. They used centered bell-shaped continuous function with compact support as an activation function (indeed, a suitable linear combination of sigmoidal
functions) and Anastassiou obtained estimations for the rate of convergence of these operators based on the modulus of continuity
of the function being approximated, that is, Jackson type approximations.

Since then, these NN operators have received considerable attention (see, e.g., \cite{costa13, costa13b, costa14}).  More recently,  Qian and Yu \cite{QiYu} provided a systematic account of the above mentioned bell-shaped activation function which unifies most of the previous works. 

\bigskip
Similar problems have also been studied in a fuzzy context by replacing the real-valued functions $C([a,b])$ by fuzzy number-valued functions $C([a,b],\fn)$. 
The main strategy so far (see, e.g., \cite{Liu}, \cite{HW}, \cite{HuWu}) has been 
 to approximate the real-valued functions $\phi_j(x)$ mentioned above by real-valued neural networks 
 taking advantage of some classical results such as \cite{cybenko, hornik}. However, this procedure adds a new layer to our network increasing the computational complexity.

Another approach is to find coefficients $\phi_j(x)$ which do not increase the computational cost, but maintaining the accuracy of the approximation.
For instance, in \cite{foma1}, a family of operators $T_n:C([0,1],\fn_{\infty})\rightarrow C([0,1],\fn_{\infty})$ defined as $$T_n(f,x):=\sum_{j=0}^n\phi_j(x)f\left(\frac{j}{n}\right),$$
where $\phi_j(x)$ are trapezoidal functions, $\sum_{j=0}^n\phi_j(x)=1$ and, for each $x\in[0,1]$, only two consecutive terms in the sum are nonzero, satisfy
$$D(T_n(f),f)=\sup_{x\in[0,1]}d_{\infty}(T_n(f,x),f(x))\leq 3{\omega}_{d_{\infty}}^{\mathcal{F}}\left(f,\frac{1}{n}\right),
$$
where ${\omega}_{d_{\infty}}^{\mathcal{F}}(f,\delta)$ stands for the fuzzy modulus of continuity (defined below).

\medskip
U. Kadak (\cite{kadak}) has recently extended Cardaliaguet-Euvrard NN operators to the context of multivariate $n$-cell fuzzy number-valued functions which are continuous with respect to the supremum metric, $d_{\infty}$. 
In this paper, we will also deal with these operators but defined between weaker types of continuous functions, namely, endograph, sendograph and level-continuous functions. 
Indeed, most studies within the topic of fuzzy-valued functions center on $d_{\infty}$-continuity, but convergence in the supremum metric is the strongest among the usual metrics in the fuzzy number set, such as the sendograph metric, $d_p$-metrics, the endograph metric or the level convergence topology. Therefore, as mentioned in \cite{HuWu}, the supremum metric continuity is a very strong condition to be placed on a fuzzy number-valued function. Maybe that is the reason for the plethora of results with this kind of functions in the literature and the lack of results when we consider other weak types of continuity.

The paper is organized as follows: Section 2 presents the main topics in the fuzzy number set $\fn$ along with the key results on sigmoidal functions and fuzzy modulus of continuity.
Section 3 is devoted to the study of fuzzy NN interpolation operators for  level continuous  functions. Section 4 deals with  sendograph metric continuous functions and Section 5 focuses on endograph metric continuous functions. We would like to remark that most properties concerning endograph metric continuity (including an example of an endograph metric continuous function which is not sendograph metric continuous)  do not seem to have made its way into the literature yet. Section 6 shows some numerical experiments to illustrate the theoretical results. The last section consists of a conclusion.


\section{Preliminaries}

Throughout the paper, a \textit{space} means a Tychonoff topological space. A fuzzy set on a space $X$ is a map $u:X\rightarrow [0,1]$.
For each $\lambda\in [0,1]$, the $\lambda$-level of a fuzzy set $u$ is defined by
\begin{align*}
[u]^{\lambda} & =\{x\in X\ :\  u(x)\geq \lambda\} \quad\text{if $\lambda>0$};\\
[u]^{0}&={\rm cl}_{X}\{x\in X\ :\ u(x)>0\}
\end{align*}
where ${\rm cl}_{X}(A)$ stands for the closure in $X$ of the set $A$.
A fuzzy number $u$ will be a fuzzy set on $\R$  such that

\begin{itemize}
\item[(FN1)] There exists $x\in \R$ such that $u(x)=1$;
    \item[(FN2)] $u$ is lower semicontinuous;
    \item[(FN3)] $u$ is convex: for each $\alpha\in [0,1]$, 
    $u(\alpha x+(1-\alpha)y\geq \min\{u(x),u(y)\};$
    \item[(FN4)] the support of $u$, $[u]^{0}$, is a compact set.
\end{itemize}
We denote by $\fn$ the set of all fuzzy numbers on $\R$. By the properties above, it turns out
that each $\lambda$-level, $[u]^{\lambda}$, is a closed interval in $\R$. We will denote by $u^-(\lambda)$ and $u^+(\lambda)$ the lower and upper end-points of that interval; that is, 
$$[u]^{\lambda}=[u^-(\lambda),u^+(\lambda)].$$
The celebrated theorem of Goetschel and Voxman  allows us to construct a fuzzy number by knowing its $\lambda$-levels.

\begin{theorem} {\rm (\cite{GV})} \label{th:GV}
 Let $u\in\fn$ and $[u]\sp{\lambda}=\mathopen [u\sp{-}(\lambda),
u\sp{+}(\lambda)\mathclose]$, $\lambda\in[0,1]$.
Then the pair of functions $u\sp{-}(\lambda)$ and
$u\sp{+}(\lambda)$ has the following properties:

\begin{enumerate}[(i)]
 \item $u\sp{-}(\lambda)$ is a bounded left continuous nondecreasing function
on $\mathopen (0,1\mathclose]$;

\item  $u\sp{+}(\lambda)$ is a bounded left continuous nonincreasing function
on $\mathopen (0,1\mathclose]$;

\item $u\sp{-}(\lambda)$ and $u\sp{+}(\lambda)$ are right continuous at
$\lambda=0$;

\item $u\sp{-}(1)\leq u\sp{+}(1)$.
\end{enumerate}

\noindent  Conversely, if a pair of functions $\alpha (\lambda)$ and
$\beta(\lambda)$ satisfy the above conditions (i)-(iv), then there
exists a unique $u\in\fn$ such that
$[u]\sp{\lambda}=\mathopen[\alpha(\lambda),\beta(\lambda)\mathclose]$
for each $\lambda\in [0,1]$.
\end{theorem}

The space $\R$ is embedded in $\fn$ since we can identify each real number $x\in\R$ with the fuzzy number given by the characteristic function $\chi_{\{x\}}$. 

We will denote by $\fn_{\infty}$ the metric space $(\fn,d_{\infty})$, where the supremum metric $d_{\infty}$ is given by
$$d_{\infty}(u,v) =\sup_{\lambda\in [0,1]}\max\{|u^-(\lambda)-v^-(\lambda)|,|u^+(\lambda)-v^+(\lambda)|\}.$$

It is worth noting that it can be rewritten as
$$d_{\infty}(u,v)=\sup_{\lambda\in[0,1]} d_H([u]^{\lambda},[v]^{\lambda}), $$
where $d_H$ is 
the Hausdorff distance. 

In general, the Hausdorff distance 
between nonempty closed sets of a metric space $(X,d)$ is defined by
$$d_H(A,B)=\max\{H^{\ast}(A,B),H^{\ast}(B,A)\}$$
where $H^{\ast}(A,B)=\sup_{a\in A}d(a,B)=\sup_{a\in A}\inf_{b\in B} d(a,b)$.

The Hausdorff distance can be infinite, but it is finite if $A$ and $B$  are compact sets or, at least, closed and bounded sets.

\begin{definition}
    Let $\sigma:\R\rightarrow\R$ be a sigmoidal function; that is, a function satisfying
    $$\lim_{x\to -\infty}\sigma(x)=0\quad\text{and}\quad 
    \lim_{x\to +\infty}\sigma(x)=1.$$
Fix $m\in\R^+$. According to \cite{QiYu},  we say that $\sigma\in\mathcal{A}(m)$ if
\begin{enumerate}
    \item $\sigma(x)$ is nondecreasing,
    \item $\sigma(x)=0$ for $x\leq -m$ and $\sigma(x)=1$ for $x\geq m$.
\end{enumerate}
Let $\varphi(x):=\sigma(x+m)-\sigma(x-m)$.
\end{definition}

 In particular, the well-known ramp functions, the smooth ramp sigmoidal function,  and the univariate central B-spline belong to the class $\mathcal{A}(m)$. 
The properties of $\varphi(x)$ are summarized in the following proposition:
\begin{proposition}{\rm ({\cite[Lemma 1]{QiYu}})}\label{prop:varphi}
If $\sigma\in\mathcal{A}(m)$, then
\begin{enumerate}[(i)]
\item $\varphi(x)$ is nonnegative;
\item $\varphi(x)$ is nondecreasing for $x<0$ and nonincreasing for $x>0$;
\item $supp(\varphi(x))\subseteq [-2m,2m]$;
\item $\varphi(x)+\varphi(x-2m)=1$ for $x\in [0,2m]$.
\end{enumerate}

\end{proposition}

Now, inspired by \cite{QiYu}, we can define the following fuzzy neural network operator.
\begin{definition}
    Let $f: [a,b]\rightarrow\fn$ and $\sigma\in\mathcal{A}(m)$. Then we can define, for $x\in [a,b]$ and $n\in\N$,
$$S_{n,\sigma}(f,x):=\sum_{k=0}^{n} \varphi\left( \frac{2m}{h}(x-x_k)\right) f(x_k),$$
where $x_k=a+kh$, $k=0,1,2,\ldots, n$ and $h=\dfrac{b-a}{n}$.
\end{definition}
We recall that $\R$ is embedded canonically in $\fn$ by  $x \mapsto \widetilde{x}$, where $x\in\R$ and $\widetilde{x}:=\chi_{\{x\}}\in\fn$. Therefore, if $f:[a,b]\rightarrow\R$,  we denote by $\widetilde{f}:[a,b]\rightarrow\fn$ the function defined by $\widetilde{f}(x):=\widetilde{f(x)}$. The following proposition shows that $S_n(f)$ matches the values of $f$ at each node $x_k$ and $S_n(f)$ is normalised.
\begin{proposition}
Let $\sigma\in\mathcal{A}(m)$ be a sigmoidal function and $f:[a,b]\rightarrow\fn$. Let $S_{n,\sigma}(f,x)$ be defined as above. Then,
 \begin{enumerate}[(i)]
     \item $S_{n,\sigma}(f,x_k)=f(x_k)$, $k=0,1,\ldots, n$.
\item $S_{n,\sigma}(\widetilde{1},x)=\widetilde{1}$, for all $x\in[a,b]$.
 \end{enumerate}
\end{proposition}
\begin{proof}

\begin{enumerate}[(i)]
    \item Fix $k\in\{1,\ldots,n-1\}$. Then, $\frac{2m}{h}(x_k-x_j)\geq 2m$ for $j=0,\ldots, k-1$ and 
    $\frac{2m}{h}(x_k-x_j)\leq -2m $ for $j=k+1,\ldots, n$. Therefore, by property (iii) of Proposition~\ref{prop:varphi}, we get
\begin{align*}
S_{n,\sigma}(f,x_k)&=     \sum_{j=0}^{n} \varphi\left( \frac{2m}{h}(x_k-x_j)\right) f(x_j)\\[1ex]
&=\sum_{j=0}^{k-1} \varphi\left( \frac{2m}{h}(x_k-x_j)\right) f(x_j)+
f(x_k)\varphi(0)+\sum_{j=k+1}^{n} \varphi\left( \frac{2m}{h}(x_k-x_j)\right) f(x_j)\\[1ex]
&=f(x_k)\varphi(0)=f(x_k)(\sigma(m)-\sigma(-m))=f(x_k).
\end{align*}
The cases when $k=0$ or $k=n$ run in a similar way.
\item Take any $x\in[a,b]$. 
If $x=x_k$ for some $k\in\{0,1,\ldots,n\}$, the result follows from (i). In other case, there exists $k\in\{0,1,\ldots,n-1\}$ such that
$x_k<x<x_{k+1}$.
Since  $x_{k+1}=x_k+h$, we have
$$ \frac{2m}{h}(x-x_{k+1})= \frac{2m}{h}(x-x_{k}-h)= \frac{2m}{h}(x-x_{k})-2m.$$
Now, $\frac{2m}{h}(x-x_{k})\in[0,2m]$ and, by applying  (\textit{iv}) of Proposition~\ref{prop:varphi}, we get 
$$\varphi\left( \frac{2m}{h}(x-x_k)\right)+\varphi\left( \frac{2m}{h}(x-x_{k+1})\right)=1.$$
On the other hand, 
$\frac{2m}{h}(x-x_j)\geq 2m$ for $j=0,\ldots, k-1$ and 
    $\frac{2m}{h}(x-x_j)\leq -2m $ for $j=k+2,\ldots, n$.
Therefore,
\begin{align*}
\sum_{j=0}^{n} \varphi\left( \frac{2m}{h}(x-x_j)\right)&=\sum_{j=0}^{k-1} \varphi\left( \frac{2m}{h}(x-x_j)\right)+ \varphi\left( \frac{2m}{h}(x-x_k)\right)\\[1ex]
&+\varphi\left( \frac{2m}{h}(x-x_{k+1})\right)+\sum_{j=k+2}^{n}\varphi\left( \frac{2m}{h}(x-x_j)\right)\\[1ex]
&=\varphi\left( \frac{2m}{h}(x-x_k)\right)+\varphi\left( \frac{2m}{h}(x-x_{k+1})\right)=1,
\end{align*}
which provides
\begin{align*} 
S_{n,\sigma}(\widetilde{1},x)&=
\left(\sum_{j=0}^{n} \varphi\left( \frac{2m}{h}(x-x_j)\right) \right)\widetilde{1}=\widetilde{1}.\\[1ex]
\end{align*}
\end{enumerate}
\end{proof}

 Let us recall the definition of the modulus of continuity of a  function
$f:[a,b]\rightarrow \R$:
\[
\omega(f,\delta):=\sup\{|f(x)-f(y)|\ :\ x,y\in[a,b],\ |x-y|<\delta\}
\]
and, in a similar way, we can define the fuzzy modulus of continuity of a fuzzy number-valued function 
$f:[a,b]\rightarrow(\fn,d)$, where
$d$ is any metric on $\fn$:
\[
\omega^{\mathcal{F}}_{d}(f,\delta):=\sup\{d(f(x),f(y))\ :\ x,y\in[a,b],\ |x-y|<\delta\}.
\]
Since $[a,b]$ is compact, if $f$ is continuous, then it is uniformly continuous (see, e.g., \cite{daners}) and, consequently, $\lim_{\delta\to 0}\omega^{\mathcal{F}}_{d}(f,\delta)=0$.
    

\section{Level continuous mappings}

U. Kadak \cite{kadak} has shown that the neural network operators $S_{n,\sigma}(f,x)$ are $d_{\infty}$-close to $d_{\infty}$-continuous functions $f:[a,b]\rightarrow\fn_{\infty}$. Namely,  $$D(S_{n,\sigma}(f),f)\leq {\omega}_{d_{\infty}}^{\mathcal{F}}\left(f,\frac{1}{n}\right).$$

In this section we consider $\fn$ endowed with the level convergence topology. As mentioned in the introduction, the class of level continuous functions is strictly wider than the class of $d_{\infty}$-continuous functions (see~\cite[Example 5.1]{FH:04}), so we need new techniques to cover this kind of functions.

Let us recall that a net of fuzzy numbers $\{u_{\alpha}\}_{{\alpha}\in I}$ is said that {\em levelly converges} to $u\in\fn$ if
$$\lim_{\alpha} u_{\alpha}^-(\lambda)=u^-(\lambda)\quad \text{and}\quad \lim_{\alpha} u_{\alpha}^+(\lambda)=u^+(\lambda)$$
for all $\lambda\in[0,1]$. This convergence defines a topology on $\fn$ with the following local base of neighbourhoods:
$$V(u,\{\lambda_1,\ldots,\lambda_n\},\vr)=
\{v\in\fn\ :\ \max_{1\leq j\leq n}\{|u^-(\lambda_j)-v^-(\lambda_j)|,|u^+(\lambda_j)-v^+(\lambda_j)|\}<\vr\},$$
where $u\in\fn$, $\{\lambda_1,\ldots,\lambda_n\}\subset [0,1]$ and $\vr>0$.

Let $\fn_{\ell}$ denote the set of fuzzy numbers endowed with the level convergence topology.

The level continuity of $f:[a,b]\rightarrow\fn_{\ell}$ means that, for every $\lambda\in[0,1]$, the real valued functions
$$f(\cdot)^{-}(\lambda):[a,b]\rightarrow\R\quad\text{and}\quad f(\cdot)^{+}(\lambda):[a,b]\rightarrow\R,$$
defined by $x\mapsto f(x)^{-}(\lambda)$ and $x\mapsto f(x)^{+}(\lambda)$, respectively, are continuous. We say that $f$ is level uniformly continuous if these two functions are uniformly continuous for every $\lambda\in [0,1]$. That is, for every $\lambda_0\in [0,1]$ and  $\vr>0$, there exists $\delta>0$ such that 
\begin{align}
|f(x)^{-}(\lambda_0)-f(y)^{-}(\lambda_0)|&<\vr, \quad\text{if $|x-y|<\delta$,}\label{eq:luc1}\\
|f(x)^{+}(\lambda_0)-f(y)^{+}(\lambda_0)|&<\vr, \quad\text{if $|x-y|<\delta$.}
\label{eq:luc2}
\end{align}
As in the metric case, it is straightforward that the continuity of a function $f:[a,b]\rightarrow\fn_{\ell}$ is equivalent to the uniform continuity as well.
\begin{proposition}
Let $f:[a,b]\rightarrow\fn_{\ell}$ be a level continuous function. Then, $f$ is level uniformly continuous.    
\end{proposition}

So, we can define a levelwise fuzzy modulus of continuity for functions $f:[a,b]\rightarrow\fn_{\ell}$ as follows:
for every $\lambda\in[0,1]$, define
$$\omega^{\mathcal{F}}(f,\delta,\lambda):=\sup\{d_H\left([f(x)]^{\lambda},[f(y)]^{\lambda}\right)\ :\ x,y\in[a,b]; \ |x-y|<\delta\}$$
and we can get the following result.
\begin{proposition}
Let $f:[a,b]\rightarrow \fn_{\ell}$ be a function. Then,
$f$ is level uniformly continuous if, and only if, $\lim_{\delta\to 0} \omega^{\mathcal{F}}(f,\delta,\lambda)=0$, for each $\lambda\in [0,1]$.
\end{proposition}
\begin{proof}
Assume that $f$ is level uniformly continuous. Fix $\lambda_0\in [0,1]$. Take $\vr>0$ and find $\delta_0>0$ such that Equations~\eqref{eq:luc1} and \eqref{eq:luc2} are satisfied.
Take any $0<\eta<\delta_0$. Then,
for all $x,y\in[a,b]$ with $|x-y|<\eta<\delta_0$ we have $d_H\left([f(x)]^{\lambda_0},[f(y)]^{\lambda_0}\right)<\vr$, which yields
$$\omega^{\mathcal{F}}(f,\eta,\lambda_0)=\sup\{d_H\left([f(x)]^{\lambda_0},[f(y)]^{\lambda_0}\right)\ :\ x,y\in[a,b]; \ |x-y|<\eta\}\leq\vr,$$
providing $\lim_{\delta\to 0} \omega^{\mathcal{F}}(f,\delta,\lambda_0)=0$.

Assume now that $\lim_{\delta\to 0} \omega^{\mathcal{F}}(f,\delta,\lambda)=0$, for each $\lambda\in[0,1]$. Fix $\lambda_0\in[0,1]$ and take $\vr>0$. There exists $\delta>0$ such that, for every $0<\eta<\delta$, we have
$\omega^{\mathcal{F}}(f,\eta,\lambda_0)<\epsilon$.
If $x,y\in[a,b]$ with $|x-y|<\delta/2$, then
$$
d_H\left([f(x)]^{\lambda_0},[f(y)]^{\lambda_0}\right)\leq \omega^{\mathcal{F}}(f,\delta/2,\lambda_0)<\vr,
$$
so Equations~\eqref{eq:luc1} and \eqref{eq:luc2} are satisfied and we are done.
\end{proof}

Let $C([a,b],\fn_{\ell})$ denote the space of level continuous fuzzy number valued funcions and assume that it is endowed with the topology of uniform convergence. A local basis at $f_0\in  C([a,b],\fn_{\ell})$ is formed by the sets
$$V(f_0,\{\lambda_1,\ldots,\lambda_m\},\vr):=
\{f\in C([a,b],\fn_{\ell})\ :\ \max_{1\leq j\leq m}d_H\left([f_0(x)]^{\lambda_j},[f(x)]^{\lambda_j}\right)<\vr, \ x\in[a,b]\}$$
where $\{\lambda_1,\ldots,\lambda_m\}\subseteq [0,1]$ and $\vr>0$.

\begin{theorem}\label{th:main1Level}
Let $\sigma\in\mathcal{A}(m)$. Fix $f_0\in  C([a,b],\fn_{\ell})$. Then, given any $\lambda\in[0,1]$, we have 
$$d_H\left([f_0(x)]^{\lambda},[S_{n,\sigma}(f_0,x)]^{\lambda}\right)\leq \omega^{\mathcal{F}}(f_0,h,\lambda),$$
for all $n\in\N$ and every $x\in[a,b]$, where $h=\frac{b-a}{n}$.
\end{theorem}
\begin{proof}
Choose $n\in\N$ and $\lambda\in[0,1]$. Fix   $x_0\in[a,b]$. Take $h=\frac{b-a}{n}$ and $x_k=a+kh$, $k=0,1,\ldots,n$. If $x_0=x_k$ for some $k$, then there is nothing to prove. On the other case, find $k\in\{0,1,\ldots, n-1\}$ such that $x_{k}<x_0<x_{k+1}$. Then,
\begin{align*}
\hspace*{-1cm}d_H\left([f_0(x_0)]^{\lambda},[S_{n,\sigma}(f_0,x_0)]^{\lambda}\right)
&=\text{\tiny $d_H\left(\left[\sum_{j=0}^{n} \varphi\left( \frac{2m}{h}(x_0-x_j)\right)f_0(x_0)\right]^{\lambda},\left[\sum_{j=0}^{n} \varphi\left( \frac{2m}{h}(x_0-x_j)f_0(x_j)\right)\right]^{\lambda}\right)$}\\[2ex]
&\leq \varphi\left( \frac{2m}{h}(x_0-x_{k})\right)
d_H\left(\left[f_0(x_0)\right]^{\lambda},\left[f_0(x_{k})\right]^{\lambda}\right)\\[2ex]
&+
\varphi\left( \frac{2m}{h}(x_0-x_{k+1})\right)
d_H\left(\left[f_0(x_0)\right]^{\lambda},\left[f_0(x_{k+1})\right]^{\lambda}\right)\\[2ex]
&\leq\text{\small
$\left(\varphi\left( \frac{2m}{h}(x_0-x_{k})\right)+\varphi\left( \frac{2m}{h}(x_0-x_{k+1})\right)\right)
\cdot \omega^{\mathcal{F}}(f_0,h,\lambda)$}\\[2ex]
&=\omega^{\mathcal{F}}(f_0,h,\lambda).
\end{align*}
This finishes the proof.
\end{proof}
As a consequence, we can show that these fuzzy neural network operators can approximate any level continuous function with any degree of accuracy.

\begin{theorem}\label{th:mainLevel}
Let $\sigma\in\mathcal{A}(m)$. Fix $f_0\in  C([a,b],\fn_{\ell})$. Then, given any $V(f_0,\{\lambda_1,\ldots,\lambda_m\},\vr)$, we can find $n\in \N$ such that 
$$S_{n,\sigma}(f_0,x)\in V(f_0,\{\lambda_1,\ldots,\lambda_m\},\vr).$$
\end{theorem}
\begin{proof}
Since $f_0$ is a level-continuous  function, we know that it is (level) uniformly continuous so we have that, for every $j=1,\ldots,m$,
$
    \lim_{\delta\to 0} \omega^{\mathcal{F}}(f_0,\delta,\lambda_j)=0$.
Thus we can find $\delta_0>0$ such that
for every $0<\delta<\delta_0$ we have 
$\max\{\omega^{\mathcal{F}}(f_0,\delta,\lambda_j)\ :\ j=1,\ldots,m\}<\epsilon$.

Choose $n\in\N$ such that $h=\frac{b-a}{n}<\delta_0$ and set $x_{k}=a+k\cdot h$, $k=0,1,\ldots, n$.  

Fix $x_0\in[a,b]$. If $x_0=x_k$ for some $k=0,1,\ldots, n$ there is nothing to prove. If this is not the case, find $k\in\{0,1,\ldots, n-1\}$ such that $x_{k}<x_0<x_{k+1}$. Then, for every $j=1,2,\ldots,m$, we have
\begin{align*}
d_H\left([f_0(x_0)]^{\lambda_j},[S_{n,\sigma}(f_0,x_0)]^{\lambda_j}\right)
\leq \omega^{\mathcal{F}}(f_0,h,\lambda_j)<\vr,
\end{align*}
which yields the conclusion.
\end{proof}

\bigskip
\section{The sendograph metric case}

Given $u\in\fn$, the sendograph of $u$ is defined by 
$${\rm send}(u):=\{(x,\alpha)\in [u]^{0}\times [0,1]\ : \ u(x)\geq\alpha\}.$$
The sendograph is a compact set and, then, we can define the 
sendograph metric by 
$$D_S(u,v)=d_H({\rm send}(u),{\rm send}(v)),\quad u,v\in\fn.$$

It is known that the convergence in the sendograph metric is weaker than that  in the level convergence topology, so that the class of $D_S$-continuous fuzzy number valued functions is strictly larger than the class of level continuous functions. Indeed, an example of this fact can be found in \cite[Example 8.1]{HuWu}.

\begin{proposition}\label{prop:bounds}
The sendograph metric on $\fn$ has the following properties:
\begin{enumerate}[(i)]
\item $D_S(\alpha u, \beta u)\leq |\alpha-\beta| \max\{|y|\ :\ y\in[u]^0\}$; for all $\alpha,\beta\in\R$, $\in\fn$.
\item $D_S(\sum_{j=1}^{n} u_j,\sum_{j=1}^{n} v_j)\leq \sum_{j=1}^n D_S(u_j,v_j)$; $u_j,v_j\in\fn$, $j=1,2,\ldots,n$.
\item If $0\leq \alpha\leq 1$, then,  for all $u,v,w\in\fn$,
$$D_S(\alpha u+(1-\alpha)v,w)\leq \sqrt{2}\; \max\{D_S(u,w),D_S(v,w)\}.$$
\end{enumerate}
\end{proposition}
\begin{proof}
The proofs can be found in \cite{HuWu}.  
\end{proof}
Now we can state the following result.

\begin{theorem}\label{th:mainSendo}
Let $\sigma\in\mathcal{A}(m)$. Fix $f_0\in  C([a,b],(\fn,D_S))$. Then, for every $n\in \N$ we have 
$$\sup_{x\in[a,b]} D_S(S_{n,\sigma}(f_0,x), f_0(x))\leq\sqrt{2}\cdot\omega^{\mathcal{F}}_{D_S}\left(f_0,\frac{b-a}{n}\right).$$

\end{theorem}
\begin{proof}
Take $n\in\N$ and fix $x_0\in[a,b]$. Take $h=\frac{b-a}{n}$ and set $x_{k}=a+k\cdot h$, $k=0,1,\ldots, n$. If $x_0=x_k$ for some $k=0,1,\ldots, n$ there is nothing to prove. So, find $k\in\{0,1,\ldots,n-1\}$ such that $x_{k}<x_0<x_{k+1}$.  By Proposition~\ref{prop:bounds}(iii) and bearing in mind that $\varphi\left( \frac{2m}{h}(x_0-x_{k+1})\right)+\varphi\left( \frac{2m}{h}(x_0-x_{k+1})\right)=1$, we have  
\begin{align*}
\hspace*{-10mm}D_S(S_{n,\sigma}(f_0,x_0), f_0(x_0))
&= D_S\left(\sum_{j=0}^{n} \varphi\left( \text{\tiny $\frac{2m}{h}(x_0-x_j)$}\right)f_0(x_j),f_0(x_0)\right)\\[2ex]
&= \text{\small $D_S\left(\varphi\left( \text{\tiny $\frac{2m}{h}(x_0-x_k)$}\right)f_0(x_k)+\varphi\left( \text{\tiny $\frac{2m}{h}(x_0-x_{k+1})$}\right)f_0(x_{k+1}),f_0(x_{0})\right)$}\\[2ex]
&\leq \sqrt{2}\cdot \max\Big\{D_S\left(f_0(x_{k}),f_0(x_{0})\right),D_S\left(f_0(x_{k+1}),f_0(x_{0})\right)\Big\}\\[2ex]
&\leq \sqrt{2}\cdot \omega^{\mathcal{F}}_{D_S}(f_0,h).\\[2ex]
\end{align*}
\end{proof}

As in the case of level continuous functions, we can now deduce a density result for fuzzy number-valued $D_S$-continuous functions.

\begin{theorem}\label{th:mainSendo2}
Let $\sigma\in\mathcal{A}(m)$. Fix $f_0\in  C([a,b],(\fn,D_S))$. Then, given any $\vr>0$, we can find $n\in \N$ such that 
$$\sup_{x\in[a,b]} D_S(S_{n,\sigma}(f_0,x), f_0(x))<\vr.$$
\end{theorem}
\begin{proof}
By applying that $f_0$ is uniformly $D_S$-continuous, we get 
$$\lim_{h\to 0}\omega^{\mathcal{F}}_{D_S}(f_0,h)=0$$
so we can find $n\in\N$ such that 
$$\omega^{\mathcal{F}}_{D_S}(f_0,\frac{b-a}{n})<\vr$$
and, by applying Theorem~\ref{th:mainSendo}, we get the conclusion.
\end{proof}

\bigskip
\section{The endograph metric case}

Another metric that we can consider on $\fn$ is the endograph metric. Recall that the endograph of a fuzzy number $u\in\fn$ is defined by 
$$end(u)=\{(x,\beta)\in\R\times[0,1]\ :\ u(x)\geq \beta\}$$
which is a closed set in $\R\times[0,1]$ and the endograph metric is defined by 
$$D_E(u,v)=d_H(end(u),end(v)), \quad u,v\in\fn,$$
where $d_H$ stands for the Hausdorff distance in $(\R\times [0,1],d)$, being $d$ the Euclidean distance.

The endograph metric satisfies that $D_E\leq D_S$ and then, the convergence in this metric is weaker (indeed, strictly weaker as we show at the end of this section) than the convergence in the sendograph metric
(see~\cite[Theorem 5]{Fang} for a relationship between the convergence in these two metrics).

Since $end(u)=(\R\times\{0\})\cup send(u)$ for each $u\in \fn$, the sendograph and endograph metrics are expected to share some behaviour.
Next, we shall give some properties related to the endograph metric similar to those stated in Proposition~\ref{prop:bounds} for the sendograph metric, which, as far as we know, had not appeared in the literature yet. First we prove three technical lemmas.

\begin{lemma}\label{prop:convex}
Let $u,v$ be two fuzzy numbers. Let $w=\alpha\cdot u+(1-\alpha)\cdot v$, for $\alpha\in[0,1]$.Then, for each $\lambda\in[0,1]$, we have that $z\in[w]^{\lambda}$ if, and only if, there exist $x\in[u]^{\lambda}$ and $y\in[v]^{\lambda}$, such that  $z=\alpha x+(1-\alpha)y$.
\end{lemma}

\begin{proof}
Assume firstly that $x\in[u]^{\lambda}$, $y\in[v]^{\lambda}$ and $z=\alpha x+(1-\alpha)y$. Then, by the convexity of fuzzy numbers, we have
\begin{align*}
   w(z)&= (\alpha\cdot u+(1-\alpha)\cdot v)(z)=\alpha u(z)+(1-\alpha)v(z)\\
    &\geq \alpha \min\{u(x),u(y)\}+(1-\alpha)\min\{v(x),v(y)\}\\
    &\geq \alpha\lambda+(1-\alpha)\lambda=\lambda.
\end{align*}
So, $z\in[w]^{\lambda}$.

Conversely, fix $z\in[w]^{\lambda}=[w^-(\lambda),w^+(\lambda)]$. There exist $0\leq k\leq 1$, such that $z=kw^-(\lambda)+(1-k)w^+(\lambda)$. Take 
$x=ku^-(\lambda)+(1-k)u^+(\lambda)\in[u]^{\lambda}$ and 
$y=kv^-(\lambda)+(1-k)v^+(\lambda)\in[v]^{\lambda}$. Then,
\begin{align*}
   z&=kw^-(\lambda)+(1-k)w^+(\lambda)\\
   &=k(\alpha u^-(\lambda)+(1-\alpha)v^-(\lambda))+(1-k)(\alpha u^+(\lambda)+(1-\alpha)v^+(\lambda))\\
   &=\alpha(ku^-(\lambda)+(1-k)u^+(\lambda))+(1-\alpha)(k v^-(\lambda)+(1-k)v^+(\lambda))\\
   &=\alpha x+(1-\alpha)y.
\end{align*}
That finishes the proof.
\end{proof}


\bigskip
\begin{lemma}\label{lem:endachiv1}
Let $K$ be a compact set in $(\R\times [0,1],d)$ and set $M=(\R\times\{0\})\cup K$. Then, for each $(x,\lambda)\in \R\times [0,1]$ we have that
$$d((x,\lambda),M)=\min\{d((x,\lambda),(x,0)),d((x,\lambda),K)\}.$$
Moreover, for each $(x,\lambda)\in \R\times [0,1]$ we have that there exists $(y,\beta)\in M $ such that
$$d((x,\lambda),(y,\beta))=d((x,\lambda),M).$$
\end{lemma}
\begin{proof}
Since $d((x,\lambda),M)=\inf_{(y,\beta)\in M} d((x,\lambda),(y,\beta))$, we have $$d((x,\lambda),M)\leq \min\{d((x,\lambda),(x,0)),d((x,\lambda),K)\}.$$
Now, given $\vr>0$, there exists $(y,\beta)\in M$ with 
$$d((x,\lambda),(y,\beta))<d((x,\lambda),M)+\vr.$$
If $\beta=0$, then we have
$$\min\{d((x,\lambda),(x,0)),d((x,\lambda),K)\}\leq d((x,\lambda),(x,0))\leq d((x,\lambda),(y,0))<d((x,\lambda),M)+\vr.$$
If $\beta\neq 0$, then $(y,\beta)\in K$
and again
$$\min\{d((x,\lambda),(x,0)),d((x,\lambda),K)\}\leq d((x,\lambda),K)\leq d((x,\lambda),(y,\beta))<d((x,\lambda),M)+\vr.$$
So we get
$$\min\{d((x,\lambda),(x,0)),d((x,\lambda),K)\}<d((x,\lambda),M)+\vr,$$
for all $\vr>0$, which yields 
$$\min\{d((x,\lambda),(x,0)),d((x,\lambda),K)\}\leq d((x,\lambda),M)$$
and we are done.

Finally, since
$$d((x,\lambda),M)=\min\{d((x,\lambda),(x,0)),d((x,\lambda),K)\},$$
we have that, if $d((x,\lambda),M)=d((x,\lambda),(x,0))$, then $(y,\beta)=(x,0)$ gives the conclusion. If 
$d((x,\lambda),M)=d((x,\lambda),K)$, the compactness of $K$ gives us a point $(y,\beta)\in K$ such that
$d((x,\lambda),K)=d((x,\lambda),(y,\beta))$ and we are done.
\end{proof}

\bigskip
It is worth noting that with the notation of the lemma above, the distance $d((x,\lambda),M)$ is always achieved at some point in $M$, despite $M$ not being a compact set.

\begin{remark}\label{rem:haus}
By applying Lemma~\ref{lem:endachiv1}, we can state the following useful result: for $u,v\in \fn$, 
$H^{\ast}(end(u),end(v))=\displaystyle{\sup_{(x,\lambda)\in end(u)}\inf_{(y,\beta)\in end(v)} d((x,\lambda),(y,\beta))\leq K}$ for some $K>0$ if, and only if,  for every $(x,\lambda)\in end(u)$ there exists $(y,\beta)\in end(v)$ with $d((x,\lambda),(y,\beta))\leq K$. 
Therefore, in order to prove that $d_H(end(u),end(v))\leq K$, we need to check two conditions:
\begin{equation}\label{eqdh1}
\forall (x,\lambda)\in end(u) \quad \exists (y,\beta)\in end(v) \ :\ d((x,\lambda),(y,\beta))\leq K    
\end{equation}
and
\begin{equation}\label{eqdh2}
\forall (y,\beta)\in end(v) \quad \exists (x,\lambda)\in end(u) \ :\ d((x,\lambda),(y,\beta))\leq K.    
\end{equation}
\end{remark}


\begin{lemma}\label{lem:convexe}
Let $u,v$ be two fuzzy numbers. 
\begin{enumerate}[(i)]
\item If $(x_1,\beta_1), (x_2,\beta_2) \in end(u)$, then
$$(\alpha x_1+(1-\alpha)x_2,\min\{\beta_1,\beta_2\})\in end(u),$$
for each $\alpha\in [0,1]$.

\item Let $w=\alpha\cdot u+(1-\alpha)\cdot v$, for some $\alpha\in[0,1]$. If $(y_1,\lambda_1)\in end(u)$ and $(y_2,\lambda_2)\in end(v)$, then  $$(\alpha y_1+(1-\alpha)y_2,\min\{\lambda_1,\lambda_2\})\in end(w).$$  

\end{enumerate}
\end{lemma}
\begin{proof}
For (ii), assume that $(y_1,\lambda_1)\in end(u)$ and $(y_2,\lambda_2)\in end(v)$. We have three  possibilities:
\begin{itemize}
\item $y_1\in[u]^{\lambda_1}$ and $y_2\in[v]^{\lambda_2}$. 
By Proposition~\ref{prop:convex}, 
$\alpha y_1+(1-\alpha)y_2\in[w]^{\min\{\lambda_1,\lambda_2\}}$, which yields $$(\alpha y_1+(1-\alpha)y_2,\min\{\lambda_1,\lambda_2\})\in end(w).$$

\item $(y_1,\lambda_1),(y_2,\lambda_2)\in\R\times\{0\}$.
Then $\min\{\lambda_1,\lambda_2\}=0$ and so $$(\alpha y_1+(1-\alpha)y_2,\min\{\lambda_1,\lambda_2\})\in \R\times\{0\} \subset end(w).$$

\item $(y_1,\lambda_1)\in\R\times\{0\}$ and $y_2\in[v]^{\lambda_2}$ (or $(y_2,\lambda_2)\in\R\times\{0\}$ and $y_1\in[u]^{\lambda_1}$).
In this case $\min\{\lambda_1,\lambda_2\}=0$ and, then,  
 $(\alpha y_1+(1-\alpha)y_2,\min\{\lambda_1,\lambda_2\})\in\R\times\{0\}\subset end(w)$.  
 \end{itemize}

A similar argument shows (i).
\end{proof}

\bigskip
The following proposition sums up Lemmas~4.1~and~4.2 in \cite{HuWu} and it is the key to get the desired result.

\begin{proposition}\label{prop:distances}
    Let $d$ be the Euclidean distance in $\R^2$.
    Let $b,k_1,k_2\in\R$ and $\alpha\in[0,1]$.
\begin{enumerate}
\item For each $x,y_1,y_2\in\R$, we have
$$d((x,b),(\alpha y_1+(1-\alpha )y_2,\min\{k_1,k_2\}))\leq \sqrt{d((x,b),(y_1,k_1))^2+d((x,b),(y_2,k_2))^2}.$$

\item Let $x_1,x_2,y_1,y_2\in\R$ and let $x=\alpha x_1+(1-\alpha )x_2$. Then,
$$d((x,b),(\alpha y_1+(1-\alpha )y_2,\min\{k_1,k_2\}))\leq \sqrt{d((x_1,b),(y_1,k_1))^2+d((x_2,b),(y_2,k_2))^2}.$$

\end{enumerate}
\end{proposition}

The following result establishes some properties of the endograph metric similar to those of the sendograph metric set out in Proposition~\ref{prop:bounds}.
\begin{proposition}\label{prop:endproperties}
The endograph metric on $\fn$ enjoys the following properties.

\begin{enumerate}[(i)] 
\item Let $u$ be in $\fn$ and let $\alpha,\beta\in\R$. Then 
$$D_E(\alpha u,\beta u)\leq |\alpha-\beta|\max\{|z|\ :\ z\in[u]^{0}\}.$$

\item
Let $u_j,v_j\in \fn$, $j=1,2,\ldots, n$ and set $u=\sum_{j=0}^{n}u_j$ and $v=\sum_{j=0}^{n}v_j$. Then
$$D_E(u,v)\leq \sum_{j=0}^{n} D_E(u_j,v_j).$$

\item 
Let $u,v,w$ be in $\fn$ and $\alpha\in[0,1]$. Then,
    $$D_E(\alpha\cdot u+(1-\alpha)\cdot v,w)\leq
   \sqrt{D_E(u,w)^2+D_E(v,w)^2}. $$

\item 
     Let $u,v,w \in\fn$ and $\alpha\in[0,1]$. Then,
    $$D_E(\alpha\cdot u+(1-\alpha)\cdot v,w)\leq
   \sqrt{2}\max\{D_E(u,w),D_E(v,w)\}. $$

\end{enumerate}

\end{proposition}
\begin{proof}
(i)
Let $(x,\lambda)\in end(\alpha u)$.  
If $(x,\lambda)\in \R\times|{0|}$, then $(x,\lambda)\in end(\beta u)$ and
$$d((x,\lambda), (x,\lambda))=0\leq |\alpha-\beta|\max\{|z|\ :\ z\in[u]^{0}\}.$$

If $(x,\lambda)\in send(\alpha u)$, then $x\in[\alpha u]^{\lambda}=\alpha[u]^{\lambda}$, so there exists  $y\in[u]^{\lambda}$ such that $x=\alpha y$.
Now, $(\beta y,\lambda)$ lies in $end(\beta u)$ and we have that
\begin{align*}
d((x,\lambda), d(\beta y,\lambda))
&=d((\alpha y,\lambda), (\beta  y,\lambda))
=|\alpha-\beta||y|\\
&\leq   |\alpha-\beta|\max\{|z|\ :\ z\in[u]^{0}\}.
\end{align*}
So, in any case, given $(x,\lambda)\in end(\alpha u)$, there exists $(y,\gamma)\in end(\beta u)$ 
with 
$$d((x,\lambda), d(y,\gamma))\leq 
|\alpha-\beta|\max\{|z|\ :\ z\in[u]^{0}\}$$
and Equation~\eqref{eqdh1} is satisfied. 
By reversing the roles of $\alpha $ and $\beta$, we get that Equation~\eqref{eqdh2} also holds, giving the desired result. 

\bigskip
(ii)
We need to prove that the corresponding Equations \eqref{eqdh1} and \eqref{eqdh2} are satisfied. That is,
\begin{itemize}
\item Given $(x,\lambda)\in end(u)$, there exists $(y,\beta)\in end(v)$ such that $$d((x,\lambda),(y,\beta))\leq \sum_{j=0}^{n} D_E(u_j,v_j).$$   
\item Given $(y,\beta)\in end(v)$, there exists $(x,\lambda)\in end(u)$ such that $$d((x,\lambda),(y,\beta))\leq \sum_{j=0}^{n} D_E(u_j,v_j).$$ 
\end{itemize}

Let $(x,\lambda)\in end(u)$.
If $(x, \lambda)\in \R\times\{0\}$, then, $(x,\lambda)\in end(v)$ and 
$$d((x,\lambda),(x,\lambda))=0\leq \sum_{j=0}^{n} D_E(u_j,v_j).$$

If $(x,\lambda)\in send(u)$, then $x\in[u]^{\lambda} =\sum_{j=0}^n [u_j]^{\lambda}$, so there exist $x_1,x_2,\ldots,x_n\in\R$ such that 
$(x_j,\lambda)\in end(u_j)$, $j=1,2,\ldots, n$ and $x=\sum_{j=1}^n x_j$.
Now, by Lemma~\ref{lem:endachiv1}, we can find $(y_j,\beta_j)\in end(v_j)$ such that $d((x_j,\lambda),(y_j,\beta_j))=d((x_j,\lambda),end(v_j))$, $j=1,\ldots, n$. Consider $\beta=\min\{\beta_j\ :\ 1\leq j\leq n\}$ and $y=\sum_{j=0}^n y_j$. It is apparent that $(y,\beta)\in end(v)$.
Then,
\begin{align*}
d((x,\lambda),(y,\beta))&=d\left(\left(\sum_{j=0}^n x_j,\lambda\right),\left(\sum_{j=0}^n y_j,\beta\right)\right)\\[2ex]
&\leq \sum_{j=0}^n d((x_j,\lambda),(y_j,\beta_j))= \sum_{j=0}^n d((x_j,\lambda),end(v_j))\\[2ex]
&\leq \sum_{j=0}^n H^{\ast} (end(u_j),end(v_j))\leq \sum_{j=0}^n  d_H(end(u_j),end(v_j))\\[2ex]
&=\sum_{j=0}^n D_E(u_j,v_j).
\end{align*}

Reversing the roles of $u$ and $v$ we also get the second condition. This finishes the proof.

\bigskip
(iii)  Firstly, we need to check the condition~\eqref{eqdh1}. Let $c=\alpha\cdot u+(1-\alpha)\cdot v$ and take $(x,\beta)\in end(w)=\R\times\{0\}\cup send(w)$.

If $(x,\beta)\in \R\times\{0\}$, then $(x,\beta)\in end(c)$ and 
$$d((x,\beta),(x,\beta))=0\leq \sqrt{D_E(u,w)^2+D_E(v,w)^2}.$$

So, assume that  $(x,\beta)\in send(w)$. By Lemma~\ref{lem:endachiv1}, there exist $(y_1,\lambda_1)\in end(u)$ and $(y_2,\lambda_2)\in end(v)$ with 
\begin{align*}
    d((x,\beta),(y_1,\lambda_1))&=d((x,\beta), end(u))\leq D_E(u,w),\\
d((x,\beta),(y_2,\lambda_2))&=d((x,\beta), end(v))\leq D_E(v,w).
\end{align*}

By Lemma~\ref{lem:convexe}, 
 $(\alpha y_1+(1-\alpha)y_2,\min\{\lambda_1,\lambda_2\})\in end(c)$ and Proposition~\ref{prop:distances} yields 
\begin{align*}
d((x,\beta),(\alpha y_1+(1-\alpha)y_2,\min\{\lambda_1,\lambda_2\}))&\leq \sqrt{d((x,\beta),(y_1,\lambda_1))^2+d((x,\beta),(y_2,\lambda_2))^2}\\
&\leq \sqrt{D_E(u,w)^2+D_E(v,w)^2}.
\end{align*}

Therefore, the condition~\eqref{eqdh1} holds.

Now, we will check the condition~\eqref{eqdh2}. For, take $(y,\lambda)\in end(c)$. 
If $(y,\lambda)\in\R\times\{0\}$, then $(y,\lambda)\in end(w)$ as well  and 
$d((y,\lambda), (y,\lambda))=0$.

If $(y,\lambda)\in send(c)$ then, by Proposition~\ref{prop:convex}, there exist $(y_1,\lambda)\in send(u)$ and $(y_2,\lambda)\in send(v)$, such that $y=\alpha y_1+(1-\alpha)y_2$.
By Lemma~\ref{lem:endachiv1}, we can find $(x_1,\beta_1)\in end(w)$ and 
$(x_2,\beta_2)\in end(w)$ such that 
\begin{align*}
 d((y_1,\lambda),(x_1,\beta_1))&\leq D_E(u,w),\\  
 d((y_2,\lambda),(x_2,\beta_2))&\leq D_E(v,w).
\end{align*}
 
By Lemma~\ref{lem:convexe}(i), we have that  $(\alpha x_1+(1-\alpha)x_2,\min\{\beta_1,\beta_2\})\in end(w)$ 
and, by applying Proposition~\ref{prop:distances}(ii), we obtain
\begin{align*}
    d((y,\lambda),(\alpha x_1+(1-\alpha)x_2, \min\{\beta_1,\beta_2\})&\leq
    \sqrt{d((y_1,\lambda_1),(x_1,\beta_1))^2+d((y_2,\lambda_2),(x_2,\beta_2))^2}\\
    &<\sqrt{D_E(u,w)^2+D_E(v,w)^2}
\end{align*}
and, being $(y,\lambda)$ an arbitray point in $end(c)$, we get that The condition~\eqref{eqdh2} holds.
Therefore, 
$$D_E(\alpha u+(1-\alpha)v, w)=d_H(end(c),end(w))\leq \sqrt{D_E(u,w)^2+D_E(v,w)^2}$$
and the proof is done.

\bigskip
(iv) An straightforward consequence of (iii).

\end{proof}

The bound in Proposition~\ref{prop:endproperties}(iv) allows us to state the main theorem in this section as follows.




\begin{theorem}\label{th:mainEndo}
Let $\sigma\in\mathcal{A}(m)$. Fix $f_0\in  C([a,b],(\fn,D_E))$. Then, for every $n\in \N$ we have 
$$\sup_{x\in[a,b]} D_E(S_{n,\sigma}(f_0,x), f_0(x))\leq\sqrt{2}\cdot\omega^{\mathcal{F}}_{D_E}\left(f_0,\frac{b-a}{n}\right).$$

\end{theorem}
\begin{proof}
The proof mimics the one given for Theorem~\ref{th:mainSendo} thanks to the above proposition.    
\end{proof}


\bigskip
We close this section with 
an example of an endograph metric continuous function which is not sendograph metric continuous. This example  confirms that the class of
endograph metric continuous functions strictly contains the class of sendograph metric continuous functions.

\begin{example}\label{ex:endNOTsend}
Let $f:[0,1]\rightarrow \fn$ be defined as
$f(t)=u_t$, $t\in[0,1]$ where
$$u_t(x)=
\begin{cases} 1 & \text{if $x\in[t,1]$,}\\
0 & \text{otherwise}.
\end{cases}$$
for $t\in(0,1]$ and 
$$u_0(x)=
\begin{cases} 1 & \text{if $x=0$,}\\
0 & \text{otherwise}.
\end{cases}$$
We have that
$[u_t]^{\lambda}=[t,1]$
for all $\lambda\in[0,1]$ and each $t\in(0,1]$ and $[u_0]^{\lambda}=\{0\}$ for all $\lambda\in [0,1]$.

\textsc{Claim:} For all $s,t\in[0,1]$, $D_E(u_t,u_s)=|t-s|$.

We can assume, without loss of generality, that $0\leq t<s\leq 1$. Firstly, we need to check that the conditions \eqref{eqdh1} and \eqref{eqdh2} are satisfied.
\begin{itemize}
\item Firstly, we will check the condition~\eqref{eqdh1}. Take $(x,\lambda)\in end(u_t)$. We need to find $(y,\beta)\in end(u_s)$ with 
$$d((x,\lambda),(y,\beta))\leq |t-s|=s-t.$$

If $\lambda=0$ there is nothing to show, since $(x,0)\in\R\times\{0\}\subset end(u_s)$. If  $\lambda>0$, then $(x,\lambda)\in send(u_t)$ which means $x\in[u_t]^{\lambda}$. In this case we distinguish between two possibilities:
\begin{itemize}
    \item[(i)] if $t=0$, then $x\in[u_0]^{\lambda}=\{0\}$ and, since $s\in[u_s]^{\lambda}=[s,1]$, we can take $(s,\lambda)\in end(u_s)$ satisfying $$d((0,\lambda),(s,\lambda))=s=s-t.$$
    \item[(ii)] if $t\neq 0$, then $x\in[u_t]^{\lambda}=[t,1]$ and $x$ could be located in $[s,1]=[u_s]^{\lambda}$ or $x\in [t,s)$. In the first case  $(x,\lambda)\in end(u_s)$ and in the second case we can take again $(s,\lambda)\in end(u_s)$ yielding
    $$d((x,\lambda),(s,\lambda))=s-x\leq s-t.$$
   
\end{itemize}
Therefore, we get $H^{\ast}(end(u_t),end(u_s))\leq s-t$.
\item Now we will check the condition~\eqref{eqdh2}. Take now $(x,\lambda)\in end(u_s)$. 
As usual we only need to consider the case $\lambda>0$. Since $[u_s]^{\lambda}=[s,1]\subset [t,1]=[u_t]^{\lambda}$, we get 
$H^*(end(u_s),end(u_t))=0$.
\end{itemize}

As a consequence, we have shown that $D_E(u_t,u_s)\leq |s-t|$.

\medskip
 On the other hand, assuming $0\leq t<s\leq 1$, we have $t\in [u_t]^{s-t}=[t,1]$ so that $(t,s-t)\in  end(u_t)$ and 
\begin{align*}
    D_E(u_t,u_s)&=d_H(end(u_t),end(u_s))\geq 
\sup_{(x,\lambda)\in end (u_t)} d((x,\lambda),end(u_s))\\
&\geq d((t,s-t),end(u_s))=\min\{s-t,d((t,s-t),send(u_s))\}=s-t,\\
\end{align*}
since 
\begin{align*}
 d((t,s-t),send(u_s))&=
\inf_{(y,\beta)\in send(u_s)} d((t,s-t),(y,\beta))\\
&=\inf_{y\in[s,1],\,\beta\in[0,1]} d((t,s-t),(y,\beta))\\
&=\inf_{y\in[s,1]\,\beta\in[0,1]} \sqrt{|t-y|^2+|s-t-\beta|^2}\\
&\geq \inf_{y\in[s,1]} |t-y|=s-t.
\end{align*}

Therefore, $D_E(u_t,u_s)\geq|t-s|$ and the claim is proved.

\medskip
The claim shows that $f$ is $D_E$-continuous.
We will show that $f$ is not $D_S$-continuous.

\medskip
Consider the sequence $t_n=\frac{1}{n}\in[0,1]$,   and the sequence $u_n=f(t_n)=u_{t_n}\in \fn$, $n=1,2,\ldots$.
Since $f$ is $D_E$-continuous we know that
$\lim_{n} D_E(u_n,u_0)=0$. But,
$$d_H([u_n]^{0},[u_0]^0)=d_H([1/n,1],\{0\})=1.$$
Therefore, the sequence $u_n=f(t_n)$ does not converge to $u_0=f(0)$ under the metric $D_S$ (\cite[Theorem 5]{Fang}) and $f$ cannot be a $D_S$-continuous function.

\end{example}




\section{Examples}

In this section, we give some numerical experiments to
illustrate the theoretical results. Namely, we study the behaviour of Cardaliaguet-Euvrard NN operators when we consider a level continuous function (borrowed from \cite{FH:04}) which is not $d_{\infty}$-continuous.

\begin{example}\label{ex:levelcont1}
Let $f:[0,1]\rightarrow \fn$ be a function defined levelwise by
\begin{align*}
f(t)^-(\lambda)&=\begin{cases}
0& \text{if $0\leq \lambda \leq \frac{1}{2}$}\\[2ex]
\left(\lambda-\frac{1}{2}\right)^t  & \text{if $\frac{1}{2} < \lambda  \leq 1,$}
\\[2ex]
\end{cases} \quad t\in(0,1];
\\[2ex]
f(0)^-(\lambda)&=\begin{cases}
0 & \text {if $0\leq \lambda \leq \frac{1}{2}$}\\[2ex]
1  & \text{if $ \frac{1}{2}<\lambda\leq 1;$}\\[2ex]
\end{cases}\\[2ex]
f(t)^+(\lambda)&=1, \quad t\in[0,1].
\end{align*}

\end{example}
We list below the maximum error between the functions $S_n(f)$ and $f$ when we increase the number $n$, comparing the results for two sigmoidal functions: $\sigma_1(x)$, a ramp function between $-1$ and $1$, and $\sigma_2(x)$,the Heaviside function, that is
\[
\begin{aligned}
\sigma_1(x)&=\begin{cases}
    0 & \text{if $x\leq -1$} \\
    \frac{x+1}{2} & \text{if $-1<x<1$}\\
    1 & \text{if $x\geq 1$} 
\end{cases} &\qquad 
\sigma_2(x)&=\begin{cases}
    0 & \text{if $x\leq 0$} \\
    1 & \text{if $x>0$} 
\end{cases}
\end{aligned}
\]

It is worth noting that we can compute the levelwise modulus of continuity of $f$ since, for $0<\vr\leq \frac{1}{2}$ and $\lambda=\frac{1}{2}+\vr$, we have that
\begin{align*}
\omega^{\mathcal{F}}\left(f,1/n,\lambda\right)&=
\sup_{|x-y|<1/n} d_H\left( [f(x)]^{\lambda}, [f(y)]^{\lambda}\right)\\
&=\max\{\sup_{|x-y|<1/n}|\vr^x-\vr^y|,\sup_{0<y<1/n} (1-\vr^y)\} \\
&=\sup_{0<y<1/n} (1-\vr^y)=1-\vr^{1/n}
\end{align*}
Therefore, expect the approximation to be less accurate if $\vr$ is small; that is, when $\lambda$ is close enough to $\frac{1}{2}$.

We will show how the error behaves for some levels: $\lambda=0.50005$, $\lambda=0.6$ and $\lambda=0.8$ (it is clear that for $\lambda\leq \frac{1}{2}$ there is no error).

Since $\sup_{x\in[0,1]}d_H([S_n(f,x)]^{\lambda},[f(x)]^{\lambda})=\sup_{x\in[0,1]} |S_n(f(x)^-(\lambda))-f(x)^-(\lambda)|$, we shall compute the maximum error by using a 10000 discretization points in $[0,1]$ and run it in a Matlab routine on a 3,6 GHz Quad-Core Intel Core i7 processor.
The results are summarized in the following tables.
\begin{table}[H]
\centering
\begin{tabular}{|c|c|c|c|}
\hline
  n   & Max. error $\lambda=0.50005$ & Max. error $\lambda=0.6$& Max. error $\lambda=0.8$\\
 \hline
   $2$  & $0.477274464149295$ &$0.096642906363681$ & $0.033863181751644$ \\
     $6$  & $ 0.160771268447668$ & $ 0.015257412531789$ & $0.004557778630386$ \\
     $10$  &$0.076773442016527$ & $0.005915352176593$  & $0.001706767168147$\\
   $50$  & $  0.004446355529542$ & $2.590718534548619e-04$& $ 7.161098352004291e-05$\\
    $100$  & $0.001166995550397$ & $6.551579516633765e-05$ & $ 1.801068345708146e-05$\\
    $1000$  & $ 1.219934826079960e-05$ & $6.619747394687181e-07$ & $1.810847709560193e-07$\\
    \hline
\end{tabular}
    \caption{$\sigma_1(x)$, $m=1$}
    \label{fig:sigma1}
\end{table}

\begin{table}[H]
\centering
\begin{tabular}{|c|c|c|c|}
\hline
  n   & Max. error $\lambda=0.50005$ & Max. error $\lambda=0.6$ & Max. error $\lambda=0.8$\\
 \hline
     $10$  &$0.389961648339656$ & $ 0.108554063832944$ & $ 0.058314808017082 $ \\
   $50$  & $0.093400355140482$ & $0.022539964353849$ & $0.011849754859139 $ \\
    $100$  & $0.047373090181003$ & $0.011314735474151$ & $    0.005965161136717$ \\
    $1000$ & $  0.004914612625552$  & $0.001149191105674$ & $6.013829067220700e-04 $ \\
    \hline
\end{tabular}
\caption{$\sigma_2(x)$, $m=1$}
    \label{fig:sigma2}
\end{table}


\section{Conclusions}
We studied neural network operators of Cardaliaguet-Euvrard type for continuous functions from a compact interval into the set of fuzzy numbers endowed with the level convergence topology and with the topology induced by the sendograph and the endograph metrics. Regarding the latter, we have shown some properties of the endograph metric that, to the best of our knowledge, were not present in the literature yet. Similarly, we have provided an example of a continuous endograph function which is not a continuous sendograph function since we have not found such an example in the literature either. A rate of convergence, by means of the modulus of continuity, was achieved in all these cases.  Finally, we add the results of some numerical experiments to illustrate the behaviour of the approximations.


\begin{thebibliography}{12}

\bibitem{anas} Anastassiou, G.A., \textit{ Rate of convergence of some neural network operators to
the univariate case}. J. Math. Anal. Appl. 212, 237-262 (1997).





\bibitem{carda} Cardaliaguet, P.; Euvrard, G., \textit{ Approximation of a function and its derivative
with a neural network}. Neural Netw. 5(2), 207–220 (1992).



\bibitem{costa13}
Costarelli, D.; Spigler, R.,
\textit{
Approximation results for neural network operators activated by sigmoidal functions}, 
Neural Networks, \textbf{44},(2013) 101--106.

\bibitem{costa13b}
 Costarelli, D.; Spigler, R.,
\textit{
Multivariate neural network operators with sigmoidal activation functions}, 
Neural Networks, \textbf{48},(2013) 72--77.


\bibitem{costa14}
Costarelli, D.,
\textit{
Interpolation by neural network operators activated by ramp functions}, 
J Math. Analysis and Applications, \textbf{419},(2014) 574--582.

\bibitem{cybenko} Cybenko, G.,\textit{ Approximation by superpositions of sigmoidal functions}. Mathematics of Control, Signals, and Systems (1989).


\bibitem{daners}
Daners, D., \textit{Uniform continuity of continuous functions on compact metric spaces}, The American Mathematical Monthly, \textbf{122}, 6 (2017), pp. 592.


\bibitem{Fang}  Fang, T.;  Wang, G.,
\textit{Endographic approach on supremum and infimum of fuzzy numbers}, Information Sciences, {\bf 159} (2004),
221--231.

\bibitem{FH:04}  Fang, J.-X.; Huang, H.,
\textit{On the level convergence of a sequence of fuzzy numbers}, Fuzzy Sets and Systems, {\bf 147} 4 (2004),
417--435.

\bibitem{foma1}
 Font, J. J.; Macario, S.,
\textit{Jackson-type approximation for fuzzy-valued functions by means of trapezoidal functions}, Iranian Journal of fuzzy systems, \textbf{20}6 (2023), 49--62.

\bibitem{gal}
Gal, S.,
\textit{Approximation theory in fuzzy setting}, in Handbook of analytic computational methods in applied mathematics. Chapman\& Hall/CRC (2000), 617--666.

\bibitem{GV}  Goetschel, R.; Voxman, W.,
\textit{Elementary fuzzy calculus}, Fuzzy Sets  Sys., {\bf 18} (1986),
31-42.

\bibitem{hornik}
Hornik, K.; Stinchcombe, M.;  White, H., \textit{Multilayer Feedforward Networks are Universal Approximators}, Neural Networks, \textbf{2} (1989).



\bibitem{HW}  Huang, H.;  Wu, C.,
\textit{Approximation of level fuzzy-valued functions by multilayer regular fuzzy neural networks}, Math. Comp. Model., {\bf 49} (2009),
1311--1318.



\bibitem{HuWu}
Huang, H.;  Wu, C.,
\textit{
Approximation of fuzzy-valued functions by regular fuzzy neural
networks and the accuracy analysis},
Soft. Comput. \textbf{18} 12 (2014), 2525--2540.

\bibitem{kadak}
 Kadac, U., 
\textit{
Multivariate fuzzy neural network interpolation operators and applications to image processing},
Expert Systems With Applications \textbf{206}  (2022), 117771.

\bibitem{Liu}  Liu, P. Y.,
\textit{Universal approximations of continuous fuzzy-valued
function by multi-layer regular fuzzy neural networks}, Fuzzy Sets and Systems, {\bf 119} (2001),
313--320.



\bibitem{QiYu}
 Qian, Y.;  Yu, D.
\textit{
Rates of approximation by neural network interpolation
operators},
Applied Mathematics and Computation \textbf{418} (2022), 126781.



\end{thebibliography}
\end{document}